\documentclass{elsunart}

\usepackage{amssymb}
\usepackage{cite}

\begin{document}

\begin{frontmatter}

  
  \title{PackLib$^2$: An Integrated Library of Multi-Dimensional
    Packing Problems} 

\author{{S\'andor P. Fekete}}
\ead{s.fekete@tu-bs.de@tu-bs.de}
\and 
\author{{Jan C. van der Veen\thanksref{jan}}}
\ead{j.van-der-veen@tu-bs.de}

\address{Department of Mathematical Optimization,
Braunschweig University of Technology,
D--38106 Braunschweig,
Germany.}

\begin{abstract}
We present PackLib$^2$, the first fully integrated benchmark library
for multi-dimen\-sio\-nal packing instances. 
PackLib$^2$ combines a systematic collection of
all benchmark instances from previous literature with a 
well-organized set of new and challenging large instances. 
The XML format allows linking basic benchmark data with other important
properties, like bibliographic information, origin,
best known solutions, runtimes, etc. Transforming instances into 
a variety of existing input formats is also quite easy, 
as the XML format lends itself to easy conversion; for this purpose,
a number of parsers are provided. Thus, PackLib$^2$ aims at becoming
a one-stop location for the packing and cutting community:
In addition to fair and easy comparison of algorithmic work
and ongoing measurement of scientific progress, 
it poses numerous challenges for future research.
\end{abstract}



\thanks[jan]{Supported by Deutsche
  Forschungsgemeinschaft (DFG) for project ``ReCoNodes'', grants FE 407/8-1 and FE 407/8-2.}

\end{frontmatter}

\section{Introduction} \label{sec:Introduction}

A crucial feature of most important real-world operations research
problems is that they tend to be computationally hard: It is highly
unlikely that there exists an ``ideal'' algorithm that will construct
a provably optimal solution in relatively short time, no matter what
instance it is faced with.  Fast and easy heuristics may return
solutions that are quite poor for ``difficult'' instances; even
sophisticated methods that are guaranteed to find an optimum may
return a solution only after prohibitively long time.  Nevertheless,
such difficult problems need to be dealt with, so algorithms have to
be constructed, tested, improved and compared.

A natural approach for evaluating the practical performance of
solution methods is to run experiments on test instances. This is even
true for problems that allow a theoretically ``good'', i.e.,
polynomial, algorithm, as this does not guarantee a useful running
time in practical applications. Obviously, the choice of test
instances may have a crucial impact on the results of such
experiments, and comparisons between alternative approaches are only
possible when similar test instances are used. Finally, beyond the
performance of individual algorithmic implementations, keeping track
of scientific progress over time is of vital interest for a research
community.  This makes it also desirable to maintain a canon of open
challenge problems that can serve as catalysts for future
developments.

\section{Benchmark Libraries}
A well-established answer to the demands described above is
establishing and maintaining benchmark libraries for a large variety
of problems. In operations research, one of the first such efforts was
undertaken by Beasley with ORLIB~\cite{b-orldt-90}, a collection of
instances for 98 different classes of combinatorial optimization
problems, ranging from airport capacity allocation problems to vehicle
routing problems, and including many different cutting and packing
instances.

Arguably the most prominent of benchmark libraries for combinatorial
optimization is the TSPLIB \cite{r-tspli-91} by Reinelt. As the name
indicates, this collection is comprised of instances of the traveling
salesman problem (TSP), even though there are also some instances of
the capacitated vehicle routing problem. The TSPLIB has been extremely
successful in various ways.  Its instances have been used for a large
variety of problems, e.g., for matching \cite{cr-cmwpm-99} 
or for finding long tours \cite{fmrt-gfhlgmmmtsp-01}.
Moreover, solving a large, previously
unsolved TSPLIB instance has become a major scientific achievement,
requiring years of work by outstanding researchers and decades of CPU
time, announced in newspaper headlines and 
in one case~\cite{abcc-stsp-98} recognized as work worthy of a 
Beale-Orchard-Hays Prize for Excellence in Computational 
Mathematical Programming. 
Clearly, the TSPLIB
demonstrates that a benchmark library can be more than just a basis
for runtime comparison of algorithmic code.

Two other examples of benchmark libraries are the MIPLIB for mixed
integer programming (presented by Bixby et al.~\cite{bbdi-mipli-93})
that serves as a benchmark for all modern integer linear programming solvers,
and the SteinLib by Koch et al.~\cite{kmv-sluls-00} that clearly demonstrates
the advances in preprocessing and exact solution methods for Steiner tree
problems.

\section{Packing Problems}
Problems of cutting and packing are among the most important problems
in both mathematical programming and real-world operations research.
Even the basic one-dimensional versions of problems like bin packing
and knapsack (discussed and used as examples in any introductory
course in optimization) are NP-hard, but are more or less
well-understood by means of linear and integer programming.
Multi-dimensional generalizations face
additional difficulties, as a straightforward modeling as a compact
integer linear program is no longer available (see Fekete and Schepers
\cite{fs-cchdop-04} for a discussion.)  As demonstrated in this
special issue (and its predecessors \cite{dw-cp-90}, \cite{bw-cp-95},
and \cite{yw-cp-02}), this gives rise to a multitude of algorithmic
approaches, dealing with a variety of problem variants.  But unlike
the progress made for the TSP, the Steiner tree problem, and for
one-dimensional packing problems, solution methods for
multi-dimensional packing problems have failed to provide
breakthroughs, where the size of solved benchmark instances has grown
by several orders of magnitude, e.g., reaching the 24978 cities of Sweden
for the TSP. At
this point, the two-dimensional knapsack instance {\tt gcut13},
consisting of 32 rectangles, is beyond the reach of the best solution
methods by Fekete and Schepers~\cite{fs-eahdo-04} and Caprara and
Monaci \cite{cm-tdkp-04}.  At the same time, multi-dimensional packing
instances are among the most popular types of puzzles, sometimes even
surrounded by quite a bit of hype, e.g., Monckton's ``Eternity'' tiling
puzzle that was the subject of a $\pounds$ 1,000,000
prize contest \cite{eternity}.

Over the years, a number of benchmark instances for cutting and
packing have been presented and used in the scientific literature.
Beasley's ORLIB had a limited number of two-dimensional instances.
Wottawa's PACKLIB~\cite{w-padp-96} was a first attempt at establishing
a general benchmark library for multi-dimensional packing.
ESICUP~\cite{esicup} started to collect data sets from different
sources; unfortunately, these instances differ in file format, making
it harder than neccessary for researchers to facilitate them in their
research.

Other instances were created and used in the context of a variety of
research papers; see Section~\ref{sec:datasets} for an overview.  It
should be noted that most of the larger instances were originally
designed as test instances for other, more restricted problems like
guillotine cutting.  This indicates that even though the capability of
algorithms for solving instances of multi-dimensional packing problems
has grown only moderately compared to those for other problems, the
development of benchmark instances has not kept up.  This makes it
desirable to establish a collection of harder instances, allowing a
basis for comparison and an ongoing challenge for further progress.

\section{PackLib$^2$} \label{sec:packlib2}

As indicated above, there has been more than one attempt at
establishing a library of cutting and packing problems. Each of these
libraries had their advantages and their shortcomings.  From the
perspective of input, the strong point of the ORLIB has been its very
simple file format: the description of a packing instance is reduced
to numbers. On the other hand, this also constitutes a disadvantage,
as a correct parsing of an instance requires reading the corresponding
article; moreover, despite of its compact representation, ORLIB
encoding still contains some redundant information.  If a cost is
given for a box, it is almost always equal to its volume. Other
important information is omitted from the files. This has lead to
modifications of instances: In \cite{cm-tdkp-04} Caprara and Monaci
accidentally solve a slight modification of the gcut instances. The
gcut instances do not contain information how many boxes of a given
type can be packed. So Caprara and Monaci assumed that there was only
one box of each type.

The most sophisticated attempt at setting up a cutting and packing
library in terms of file format has been Wottawa's PACKLIB
\cite{w-padp-96}. Wottawa tried to promote a file format that was
self-describing. Its drawback was that, at that time, it required a
rather sophisticated and error-prone parser to read such a file.

As indicated by the name, PackLib$^2$ is a successor of PACKLIB. Like
PACKLIB it employs state-of-the-art technologies for representing
cutting and packing instances. Unlike all previous attempts,
PackLib$^2$ files not only capture instances but also references to
creators, references to attempts at solving the instances,
bibliographic information, and solution data.  Because PackLib$^2$ is
XML-based, the parser for our file format is based on standard
technology.  As distribution has progressed from electronic mail
\cite{b-orldt-90} over ftp-servers (possibly dressed in a web
interface) \cite{w-padp-96}, we are making full use of current
cross-referencing possibilities of websites: PackLib$^2$ is hosted at

\begin{verbatim}
      http://www.math.tu-bs.de/packlib2
\end{verbatim}

As mentioned above, the core of PackLib$^2$ is a set of XML-based
files, one for each article listed on the website. Each of these files
is subdivided into three sections:

\begin{enumerate}
\item The description section gives general information about the
  article.  This is the only mandatory section of a PackLib$^2$ XML
  file. This section basically is an extension of the BibTeX format.
\item If new problem instances or modifications of known instances are
  described in the article, they are listed in the problem section.
\item Finally the results section lists the computational results of
  the article. Whenever we are able to obtain a complete description
  of the solution, the solution itself and resulting images are
  available on PackLib$^2$.
\end{enumerate}

A detailed and up-to-date description (including future updates and
extensions) of the file format is available on the PackLib$^2$
website.
Based on these XML files, the PackLib$^2$ website is rebuilt
(half-)automatically whenever a new file is added. This automatism
ensures the integrity and comparability of the results obtained by
different researchers.

Besides listing instances and results, PackLib$^2$ also hosts cutting
and packing software. At this point, a parser for the XML files, as
well as converters to an ORLib-like format and to the old PackLib
format are available.  Furthermore a program that generates zero--one
ILP formulations based on~\cite{b-autdg-85} is available.

\section{Description of Data Sets} \label{sec:datasets}

In this section we describe the instances that are currently part of
PackLib$^2$, listed in chronological order. So far, all instances
are two-dimensional instances. Most instances presented were
originally posed as guillotine cutting stock problems. They have been
reused in other settings as well. We have classified all results using
the new typology presented in \cite{whs-itcpp-06}.
  
The oldest and smallest instance was defined by Herz in
\cite{h-rcptd-72}, presenting a recursive procedure for the
two-dimensional guillotine cutting problem. The algorithm was
implemented in PL/1 and one hand-crafted instance was solved.

In 1977 Christofides and Whitlock presented a tree-search algorithm
for the same problem \cite{cw-atdcp-77}. They tested their algorithm
on 7 instances of the two-dimensional guillotine cutting problem.
Three of these instances were described explicitly and are part of
PackLib$^2$. The test problems were randomly generated. Given the
container $A_0$ with area $\alpha_0 = L_0 W_0$, $m$ boxes with area
$\alpha_i$ were drawn uniformly at random from the interval $[0, 0.25
\alpha_0]$. Given these areas, the length $l_i$ of a box was drawn
from the interval $[0, \alpha_i]$ and then rounded up to the nearest
integer. The width of the boxes was calculated by $w_i = \lceil
\alpha_i / l_i \rceil$. The boxes were then weighted by $\upsilon_i =
r_i \alpha_i$, where $r_i$ is a uniformly distributed random number in
the range from 1 to 3.

Another set of instances is due to Bengtsson. In \cite{b-prpha-82} he
gave 10 two-dimensional bin-packing problems. For each of the
instances he generated 200 boxes with length $\lfloor 12 r + 1
\rfloor$ and width $\lfloor 8 r + 1 \rfloor$. Here $r$ is drawn from a
uniform distribution in the range (0,1). Two different containers of
width 25 and length 10 and of width 10 and height 25 were considered.

Two practical instances of the constrained two-dimensional cutting
stock problem taken from applications in the lumber industry were given
by Wang in \cite{w-tactd-83}.
In \cite{b-autdg-85}, Beasley introduced 13 randomly generated problems
of different sizes. Here the length $l_i$ of each box
was generated by sampling an integer from the uniform distribution
$[L_0/4, 3L_0/4]$. The width $w_i$ was drawn from the interval
$[W_0/4, 3W_0/4]$. $L_0$ and $W_0$ denote width and height of the
container. The value of the boxes was set to their area.
  
In \cite{b-etdng-85} Beasley introduces 12 more packing instances. The
procedure for generating instances is essentially the same as in
\cite{cw-atdcp-77}. In addition, each box may be packed more than
once. The maximal count was generated by sampling an integer from the
uniform distribution $[1, 3]$.

Hadjiconstantinou and Christofides introduced 12 new data sets for the
general, orthogonal, two-dimensional kanpsack problem. They were
generated as follows. The dimensions $l_i$ and $w_i$ of the boxes
$R_i$ are integers sampled from the uniform distributions $[1, 0.75
W_0]$ and $[0.1, 0.75W_0]$, respectively. The integer value
$\upsilon_i$ was generated by multiplying $l_iw_i$ by a real random
number drawn from a uniform distribution and rounding up the result to
the nearest integer.

PackLib$^2$ also hosts two randomly generated instances of the
two-dimensional cutting-stock problem by Tsch\"{o}ke and Holth\"{o}fer
\cite{th-npact-95}. 
Five more instances were listed explicitly in Hifi's article
\cite{h-ivbea-97}.
In \cite{fs-neagoddkp-97} five new two-dimensional knapsack instances
are defined. They were produced by a method described in
\cite{mv-estdf-98, mpv-tdbpp-00}.

For the (un)weighted constrained two-dimensional cutting stock problem
seven problems were given by Cung, Hifi, and Le Cun in
\cite{chc-ctdcs-00}. The box sizes $l_i$ and $w_i$ are chosen uniformly
at random from the intervals $[0.1L_0, 0.75L_0]$ and $[0.1W_0, 0.75W_0]$. The
weight assigned to the boxes is computed by $\upsilon_i=\lceil\rho
l_i w_i\rceil$, where $\rho=1$ for the unweighted case and
$\rho\in[0.25, 0.75]$ for the weighted case.

Finally there are the C-instances by Hopper and Turton
\cite{ht-eimhh-01}. These instances have 17 to 197 items. Three
instances were generated for each problem catogery. Width and height
of the boxes are produced randomly with a maximum aspect ratio of 7.
The problems were constructed such that the optimal solution is
known in advance. The ratio of the two dimensions of 
the container varies between 1 and 3.

Many of these instances have been used again and again by other
researchers to demonstrate the effectiveness of their algorithms. This
fact is now traceable through PackLib$^2$.

\section{Conclusions} \label{sec:conc}

There are various possible expansions of PackLib$^2$. An obvious one is to
add more instances; contributions are always welcome, especially
if they are provided with accompanying data, ideally in XML format.
Other enhancements include two- and three-dimensional problem 
variations (like including order constraints,
as suggested in \cite{fkt-mdpoc-01}) and the possibility to consider
other types of objective functions.

We are also in the process of soliciting practical problem instances
from the technical computer science community dealing with
reconfigurable computing, where
the objective is to place reconfigurable modules in two-dimensional space
(on a reconfigurable device like an FPGA) and time (as space may be re-used
when a module is no longer needed). This means that instances
are three-dimensional. See \cite{tfs-ohrt-01} for a more detailed description.

Another upcoming step will be to upgrade PackLib$^2$ to host more
algorithms.  In the very near future we plan to post the
implementation of our leading-edge packing code
\cite{fs-neagoddkp-97,fs-eahdo-04}, as well as meta-heuristic based
algorithms to tackle even larger instances.  Again, constributions are
welcome.

\section*{Acknowledgements}

We thank Andreas Ahrens, Peter Degenkolbe, and Christopher Tessars
for their help with setting up and maintaining PackLib$^2$, and
two anonymous referees for suggestions that helped to improve
the presentation of this paper.
 
\bibliographystyle{plain}
\bibliography{PackLib}

\end{document}